\documentclass[12pt,dvips]{article}
\usepackage{amsmath,amsfonts,amssymb,amsthm,graphicx,verbatim,subfig}
\usepackage[left=2cm, right=2cm, top=1.5cm, bottom=1.5cm]{geometry}
\usepackage{tikz}
\ifx\pdftexversion\undefined

\usepackage[a4paper,colorlinks,link=black,filecolor=black,citecolor=black,urlcolor=black,pdfstartview=FitH]{hyperref}
\else

\usepackage[a4paper,colorlinks,linkcolor=black,filecolor=black,citecolor=black,urlcolor=black,pdfstartview=FitH]{hyperref}
\fi

\font\sixbb=msbm6
\font\eightbb=msbm8
\font\twelvebb=msbm10 scaled 1095
\newfam\bbfam
\textfont\bbfam=\twelvebb \scriptfont\bbfam=\eightbb
                           \scriptscriptfont\bbfam=\sixbb
\def\bb{\fam\bbfam\twelvebb}

\newcommand{\Int}{{\bb Z}}

\newcommand{\FF}{{\bb F}}

\newtheorem{theorem}{\bf Theorem}
\newtheorem{claim}[theorem]{\bf Claim}

\newtheorem{corollary}[theorem]{\bf Corollary}
\newcommand{\enp}{\begin{flushright} $\Box$ \end{flushright}}
\newcommand{\beq}[0]{\begin{equation}}
\newcommand{\enq}[0]{\end{equation}}

\newcommand{\dn}{\Delta_{n-1}}

\newcommand{\rk}{{\rm rank}}

\newcommand{\supp}{{\rm supp}}

\newcommand{\lk}{\text{lk}}

\newcommand{\csy}{\text{csy}}
\newcommand{\hzero}{\widehat{0}}
\newcommand{\hone}{\widehat{1}}
\newcommand{\ovl}{\overline{L}}
\newcommand{\idd}{{\rm id}}
\newcommand{\str}{{\rm st}}

\newcommand{\ord}{\text{ord}}
\newcommand{\dist}{{\rm dist}}
\newcommand{\fix}{{\rm fix}}
\newcommand{\ffix}{{\rm Fix}}
\newcommand{\set}[1]{\{#1\}}
\newcommand{\Set}[1]{\left\{#1\right\}}

\newcommand{\remove}[1]{}

\title{Near Coverings and Cosystolic Expansion -  \\ an example of topological property testing}
\author{Irit Dinur \thanks{Department of Computer Science and Applied Mathematics, Weizmann Institute. e:mail: irit.dinur@weizmann.ac.il~.
This work is supported by an ERC-CoG grant.}
\and Roy Meshulam \thanks{Department of Mathematics,
Technion, Haifa 32000, Israel. e-mail:
meshulam@technion.ac.il~. Supported by ISF grant 326/16.}}

\begin{document}
\maketitle

\begin{abstract}
We study the stability of covers of simplicial complexes. Given a map $f:Y\to X$ that satisfies almost all of the local conditions of being a cover, is it close to being a genuine cover of $X$? Complexes $X$ for which this holds are called {\em cover-stable}. We show that this is equivalent to $X$ being a cosystolic expander with respect to non abelian coefficients. This gives a new combinatorial-topological interpretation to cosystolic expansion which is a well studied notion of high dimensional expansion.
As an example, we show that the $2$-dimensional spherical building $A_{3}(\FF_q)$ is cover-stable.

We view this work as a possibly first example of ``topological property testing'', where one is interested in studying stability of a topological notion that is naturally defined by local conditions.
\end{abstract}

\section{Introduction}
\label{s:intro}

Many central topological structures, e.g. vector bundles and covering spaces, are defined in terms of local conditions. Classification theorems for such structures are often formulated in cohomological terms.
For example, real line bundles over a compact space $X$ are classified by $H^1(X;\Int_2)$, while complex line bundles over $X$ are classified by $H^2(X;\Int)$. A natural challenge that arises is to formulate and prove
approximate (or stability) versions of such classification theorems. Roughly speaking, such results
would state that under suitable assumptions on $X$, if a structure satisfies
all but a small fraction of the local conditions, then it must be close to a structure that satisfies all of the local condidtions. We view such results as ``topological property testing''.\\

The notion of a covering space plays a key role in topology. In this paper we study the stability of this notion: Given a map that satisfies nearly all of the local requirements of being a covering map, is it close to a genuine covering map? Let us call complexes for which this holds {\em cover-stable}.
We show that a complex is cover-stable if and only if it is a cosystolic expander with respect to certain non-abelian coefficients. We further show that spherical buildings are such expanders, and hence are cover-stable.

Cosystolic expansion is a cohomological notion of expansion in simplicial complexes, that
came up independently in the study of random complexes \cite{LM06,MW09}, and in Gromov's remarkable work on the topological overlap property \cite{Gromov10}.
Its graphical one dimensional instantiation coincides with the smallest Cheeger constant of
a connected component of the graph. Here we show that its two dimensional instantiation is equivalent to cover stability. This gives a new combinatorial-topological interpretation for two-dimensional cosystolic expansion.

\paragraph{Covers and near-covers}
Let us introduce the notion of near covers by first discussing this notion for the more familiar case of graphs. We begin with $2$-covers in the introduction and move to the more general notion of $t$-cover in the body of the paper. A $2$-cover of a graph $G$ is a graph $G'$ with twice as many vertices and twice as many edges as $G$. Every vertex $v$ in $G$ is covered by two vertices $[v,0]$ and $[v,1]$, and every edge $\set{u,v}$ in $G$ is covered by two disjoint edges, either $\set{[u,0],[v,0]},\set{[u,1],[v,1]}$ or $\set{[u,0],[v,1]},\set{[u,1],[v,0]}$.
Equivalently, a $2$-cover of $G$ is a graph $G'$ together with a surjective $2$-to-$1$ homomorphism\footnote{A graph homomorphism $p:G'\to G$ is a map from $V(G)$ to $V(G')$ that sends edges to edges.} $p:G'\to G$ such that for each vertex $v'\in G'$, $p$ is a local isomorphism between the neighbors of $v'$ and the neighbors of its image $p(v')$.


How stable is this definition? suppose we have some simplicial mapping $G'\to G$ that is $2$-to-$1$ and surjective. It is a covering if the preimage of each edge of $G$ is a pair of disjoint edges. What if this holds only for nearly all edges? Does it necessarily mean that $G'$ is close to a covering? The answer is an easy yes: one can fix the lift each non-properly covered edge of $G$ without affecting anything else. So, this question is not very interesting. Nevertheless, things get more interesting once we move to two dimensions.

Let $X$ be a two-dimensional simplicial complex (a hypergraph with hyperdges of size at most $3$ that is downwards closed under containment).
A $2$-cover of $X$ is a simplicial complex $X'$ together with a surjective $2$-to-$1$ simplicial mapping $p:X'\to X$ such that the mapping from $X'$ to $X$ is a local isomorphism between the neighborhood of each vertex $v'\in X'$ and the neighborhood of its image $p(v')$.

To get a better sense of what it means to be a $2$-cover we note that $X'$ has twice as many vertices, twice as many edges, and twice as many triangles compared to $X$. Every vertex/edge in $X$ has two disjoint preimages just as in the case of covers of graphs. In addition, every triangle in $X$ must have now exactly two preimage triangles. For example, for a triangle $\set{u,v,w}\in X$, its preimages can be $\set{[u,0],[v,1],[w,0]}$ and $\set{[u,1],[v,0],[w,1]}$.
Note that this puts non-trivial requirements about the preimages of the edges $uv,uw,vw$: after having decided freely about how to cover $uv$ and $uw$, there is only one legal way to cover the third edge $vw$. \\

Next, we turn to near-covers. We consider all pairs $(f,X')$ where $X'$ is a simplicial complex and $f:X'\to X$ is a surjective $2$-to-$1$ simplicial map $f:X'\to X$. Moreover, we require that all vertices and edges of $G$ are covered properly by $f$. We measure the failure of $(f,X')$ to be a cover by counting the fraction of triangles that are not properly covered (we call this the deficientcy of the mapping, see \eqref{eq:defic} for the formal definition).
Note that the restriction to consider only mappings for which all vertices and edges are properly covered is not as restrictive or arbitrary as it might seem, because every  mapping can easily be modified to satisfy this requirement without affecting any properly-covered triangle. The interesting ``action'' happens only at the level of triangles.

The triangles of $X$ can be viewed as a test, in the property-testing sense, for $(f,X')$ being a genuine cover.
\paragraph{Triangle test:} Pick a random triangle, check that its preimage under $f$  consists of $2$ disjoint triangles. \\

We are interested in relating the failure of the test to the distance of the complex from being a true cover. The later is measured by how many edges of $X'$ need to be changed to move from $X'$ to a true cover $X''$.  The distance is measured by the natural Hamming distance on the edges of $X'$.
\remove{
Formally, let $f:X' \rightarrow X$ be an arbitrary surjective simplicial map between two simplicial complexes $X'$ and $X$. For a vertex $u'$ of $X'$ with an image $f(u')=u$, define the
\emph{local deficiency} of $f$ at $u'$ to be the number of edges in the link of $u$ that are not covered by an edge in the link of $u'$, namely
$$\mu_f(u')=|\{e \in \lk(X,u)(1): e \not\in f(\lk(X',u'))\}|.$$
The \emph{deficiency} of the map $f:Y \rightarrow X$ is given by
$$m_f(X')=\sum_{u \in X(0)}\frac{1}{|f^{-1}(u)|} \sum_{{u'} \in f^{-1}(u)} \mu_f({u'}).$$
The deficiency $m_f(X')$ is proportional to the fraction of triangles that are not propery covered\footnote{Failing to be covered can happen in several different ways, each adding a different constant to this count.}.
}

We say that a complex $X$ is {\em cover-stable} if the triangle test is a good test. Namely, if whenever the pair  $(f,X')$ has small deficiency, namely it passes the triangle test with high probability,  then it is close to a true cover $f'',X''$ (that passes the test with probability $1$).

\paragraph{Covers, cocycles, and expansion.}
For the sake of introduction let us describe the classical connection \cite{Steenrod51,Surowski84} between cocycles and covers.
Let $X$ be a two-dimensional simplicial complex and let $\phi:X(1)\to F_2$ be a labeling of the edges of $X$ with coefficients from $F_2$. $\phi$ is called a cocycle if for every $uvw\in X(2)$ we have $\phi(uv)+\phi(vw)+\phi(wu)=0$. From $\phi$ we can construct a simplicial complex $X_\phi$ as follows.
\begin{itemize}
\item Duplicate each vertex $u$ of $X$ to the pair of vertices $[u,0]$ and $[u,1]$.
\item Lift an edge $uv$ to the pair of edges $\set{[u,0],[v,0]},\set{[u,1],[v,1]}$ if $\phi(uv)=0$, and to the pair of edges  $\set{[u,0],[v,1]},\set{[u,1],[v,0]}$ if $\phi(uv)=1$.
\end{itemize}
This describes how to lift the vertices and edges. Having fixed those, the preimage of each triangle $uvw$ can either be two disjoint trinagles, or it can be a $6$-cycle (for exmaple: $u0-v0-w0-u1-v1-w1-u0$). In the former case we add the two triangles to $X_\phi$, in the latter case we have nothing to add.

Observe that $X_\phi$ is always a surjective simplicial map, and whenever a triangle  equation is satisfied by $\phi$, that triangle is properly covered by $X_\phi$. Thus, $\phi$ is a cocycle if and only if $X_\phi$ is a $2$-cover of $X$. Moreover, the deficiency of $X_\phi$ is exactly proportional to the number of triangle equations violated by $\phi$.\\

Cosystolic expansion relates the amount of triangles whose equation is violated by $\phi$ to the distance of $\phi$ from a true cocycle. A complex for which the former always bounds a constant multiple of the latter is called a cosystolic expander.
The notion of cocyle expansion (cosystolic expansion) was introduced  \cite{LM06,MW09,Gromov10} as a higher dimensional generalization of edge expansion in graphs. This notion has gained interest in recent years and for example was shown to imply the topological overlapping property. Cohomology is always specified with respect to specific coefficients, and so far most of the works focused on coefficients from $F_2$, the field with two elements. This type of cosystolic expansion dictates the cover-stability of $2$-covers.
For $t$-covers when $t>2$ we will describe in Section \ref{subs:csy} the notion of cosystolic expansion with respect to non-abelian group coefficients.

\begin{theorem}[Main, informal] $X$ is cover-stable if and only if $X$ is a cosystolic expander.
\end{theorem}
A more precise version of this theorem is stated as Theorem \ref{t:mainformal}.
Just like edge expansion gives a quantitative measure to the ``amount'' of connectivity of a graph, this theorem gives a new interpretation for two-dimensional cosystolic expansion as giving a quantitative measure for another combinatorial-topological property- that of being cover-stable.

Interestingly, even if $X$ is $d$-dimensional for $d>2$, the cover stability of $X$ is completely determined by the two-dimensional cosystolic expansion of $X$. This should not be too surprising because it is well known that the covers of $X$ are completely determined by the first cohomology.
As an example we show
\begin{theorem}[Informal, see formal version in Theorem \ref{t:sph}] The two-dimensional spherical building over any finite field is cover-stable.
\end{theorem}
\subsection*{Motivation}
\paragraph{Expansion and stable local to global phenomena}
Garland's method \cite{Garland73} is a general way to deduce global information about a complex by looking at the local views, more specifically, at the local structure of links (roughly speaking, neighborhoods) of a given complex. Originally Garland has shown a vanishing of the global cohomology by studying local link structure. This approach has been used in \cite{KKL14, EvraK16} to deduce so-called cosystolic expansion of Ramanujan complexes (this was further used to show that these are the first sparse complexes that have Gromov's topological overlapping property). A beautiful example is the trickling down theorem of Oppenheim \cite{OppGarland} that shows that if a two-dimensional complex has a connected $1$-skeleton and all of its links have good spectral expansion, then the $1$ skeleton must itself be a good spectral expander.

It is natural to wonder about a more stable version of this statement: what can be said when $99.9\%$ of the links are good spectral expanders? Perhaps one can hope that the complex is close to one that has an expanding $1$-skeleton? This turns out false. The trickling down theorem can badly fail if even very few of the links are not expanders. Think of two copies of an expanding complex joined by a single edge. The new complex can have excellent expansion in all links except two, yet the resulting complex has a $1$-skeleton that is very far from an expander.

Note however, that in this negative example the new complex is very close to a (disconnected) $2$-cover of the original complex. Could this always be the case?
Our theorem can be interpreted as showing a {\em stable version of the trickling down theorem}, for the class of complexes $Y$ for which there is a near covering mapping from the complex $Y$ to a complex $X$ that is a cosystolic expander (if $X$ happens to have spectrally expanding links then this condition would imply that 99\% of the links of $Y$ are spectrally expanding).

It is quite interesting to find a more general stable trickling down theorem.
One potential application for a stable trickling down theorem is towards a combinatorial construction of strong high dimensional expanders, in analogy to the zigzag construction of one-dimensional expanders \cite{RVW}. So far there are several known ``combinatorial'' constructions of high dimensional expanders \cite{Conlon2019, ChapmanLP19, siqi}, but none of these have links that are sufficiently expanding to apply the trickling down theorem.

\paragraph{Covers and agreement tests}
Another completely different motivation for studying stability of covers comes from agreement tests. These are certain property testing results that often underly PCP constructions.
In an agreement test one starts out with nearly matching local functions and the high dimensional expansion is used for stitching them together into one global function.
In \cite{dinur2017high} it was shown that high dimensional expanders support agreement tests. Although this gives a very strong derandomization for direct product tests in the so-called 99\% regime, no such derandomization is known for the (arguably more interesting) so-called 1\% regime. One can show \cite{DinurKamber} that the 1\% question is related to a list-agreement test, in which the local functions are replaced by lists of local functions and agreement is replaced by matching pairs of local lists. One can often reduce from 1\% agreement to 99\% list-agreement, and the end result is a near-cover of the underlying complex in the sense that we study here.

Thus, understanding which complexes are cover-stable can lead to new (and derandomized) 1\% agreement tests.  

\paragraph{Property testing and expansion}
Kaufman and Lubotzky \cite{KaufmanL14} gave a property testing interpretation to cosystolic expansion of a given complex $X$. Specifically, they showed that the test given by the coboundary operator is a good property tester (for the property of being a cocycle) iff $X$ is a cosystolic expander.
In this work we show another property testing interpretation for the cosystolic expansion of $X$. We show that $X$ is cover stable iff $X$ is a cosystolic expander iff our triangle test is a good property tester for the property of being a cover.

This aligns well with the general agenda that (high dimensional) expansion and testability go hand in hand.

It is interesting to continue to explore other topological notions that are defined by local conditions, and understand whether stability and local testability of these notions can be related to further notions of high dimensional expansion. In particular, our work only pertains to two-dimensional cosystolic expansion, and it is intriguing to find combinatorial interpretations for higher dimensional cosystolic expansion.

\section{Formal definitions and statements of results}
\remove{In this note we study the connection between the deficiency of a map, and its proximity (in a precise sense defined below) to a genuine covering map. In this note we establish a stability version of the well known classification of $G$-covering spaces of a
complex $X$, by the first cohomology set $H^1(X;G)$.
}
In this section we briefly recall some topological and combinatorial notions that play a role in our approach. We begin with general preliminary definitions in Subsection \ref{subs:prelim}. Subsection \ref{subs:nearcovers} introduces our precise notion of near covers. Subsection \ref{subs:nonab} is concerned with $1$-cohomology of a complex $X$ with non-abelian coefficients. In Subsection \ref{subs:ccm} we describe a classical construction that associates coverings with $1$-cohomology classes. In Subsection \ref{subs:csy} we recall the definition cosystolic expansion with non-abelian coefficients.
Finally, in Subsections \ref{subs:res} and \ref{subs:sph} we state our results.
\subsection{Preliminary definitions}\label{subs:prelim}
We start with some definitions.
Let $X$ be an $(n-1)$-dimensional pure simplicial complex on the vertex set $V$.
Let $X(k)$ denote the set of $k$-simplices of $X$, and let $X_{\ord}(k)$ denote the set of ordered $k$-simplices of $X$. Let $f_k(X)=|X(k)|$.
The \emph{star} and the \emph{link} of a simplex $\tau \in X$ are given by
\begin{equation*}
\begin{split}
\str(X,\tau)&=\{\sigma \in X: \sigma \cup \tau \in X\}, \\
\lk(X,\tau)&=\{\sigma \in \str(X,\tau): \sigma \cap \tau=\emptyset\}. \\
\end{split}
\end{equation*}
Define a weight function $c_X$ on the simplices of $X$ by
$$
c_X(\sigma)=\frac{|\{\tau \in X(n-1): \tau \supset \sigma\}|}{\binom{n}{|\sigma|} f_{n-1}(X)}=
\frac{f_{n-|\sigma|-1}(\lk(X,\sigma))}{\binom{n}{|\sigma|} f_{n-1}(X)}.
$$
For each $k$ the weights on $X(k)$ can be interpreted as a probability measure given by first choosing a top dimensional face $\sigma$ uniformly and then a $k$ face contained in $\sigma$. In particular note that $\sum_{\sigma \in X(k)} c_X(\sigma)=1$ for $0 \leq k \leq n-1$. Additionally, if $\alpha \in X$ and $\beta \in \lk(X,\alpha)$ then
\begin{equation*}
\label{e:prod}
c_X(\alpha)c_{\lk(X,\alpha)}(\beta)=\binom{|\alpha|+|\beta|}{|\alpha|}^{-1}c_X(\alpha \cup \beta).
\end{equation*}
In particular, if $v \in X(0)$ and $e \in \lk(X,v)(1)$ then
\begin{equation}
\label{e:prod1}
c_X(v)\cdot c_{\lk(X,v)}(e)=\frac{1}{3} \cdot c_X(v \cup e).
\end{equation}
\noindent
Let $Y$ be another simplicial complex and let $p:Y \rightarrow X$ be a surjective simplicial map.
The pair $(Y,p)$ is a \emph{covering} of $X$ if for any $u \in X(0)$ and $\tilde{u} \in p^{-1}(u)$, the induced mapping
$p: \str(Y,\tilde{u}) \rightarrow \str(X,u)$ is an isomorphism.
Consider now an arbitrary surjective simplicial map $f:Y \rightarrow X$ between two pure simplicial
complexes $Y$ and $X$. For
a vertex $\tilde{u}$ of $Y$ with an image $f(\tilde{u})=u$, let
$$D_f(\tilde{u})=\{e \in \lk(X,u)(1): e \not\in f(\lk(Y,\tilde{u}))\}.$$
Define the \emph{local deficiency} of $f$ at $\tilde{u}$ by
$$\mu_f(\tilde{u})=\sum_{e \in D_f(\tilde{u})} c_{\lk(X,u)}(e).$$
The \emph{deficiency} of the map $f:Y \rightarrow X$ is given by
\begin{equation}\label{eq:defic}
m_f(Y)=\sum_{u \in X(0)}\frac{c_X(u)}{|f^{-1}(u)|} \sum_{\tilde{u} \in f^{-1}(u)} \mu_f(\tilde{u}).
\end{equation}
The weights are actually only useful when the complex $Y$ is more than two-dimensional (in this case some triangles potentially have more weight than others). For a two-dimensional complex we can simplify the definition to an unweighted one,
\[ m_f(Y)=\frac{1}{3|X(2)|}\sum_{u \in X(0)}\frac{1}{|f^{-1}(u)|} \sum_{\tilde{u} \in f^{-1}(u)} |D_f(\tilde u)|.
\]

We view $m_f(Y)$ as a measure of the failure of $f:Y \rightarrow X$
to be a covering map. When no confusion can arise concerning the surjection $f$, we will abbreviate $D_f(\tilde{u}),\mu_f(\tilde{u})$ and $m_f(Y)$ by $D(\tilde{u}),\mu(\tilde{u})$ and $m(Y)$.

\remove{
Given two $t$-to-$1$ covering maps $(f_1,Y_1)$ and $(f_2,Y_2)$ let us define
\[ dist( (f_1,Y_1),(f_2,Y_2)) = \sum_{e\in X(1)} c_X(e) 1_{f_1^{-1}(e)\neq f_2^{-1}(e)}.
\] 

A different measure on how close $(f,Y)$ is to a covering map is the {\em proximity measure}. How much do we need to change $(f,Y)$ to make it into a genuine cover $(p,Y')$ (with zero deficiency). We will measure this distance in terms of how many edges of $Y$ need to be changed.

the edges of $Y$ to make $X$ is said to be $h$-cover-stable if for every surjective simplicial map $f:Y\to X$ one can bound the distance of a given $X'$ from a cover using its deficiency. }

\subsection{Near Covers}\label{subs:nearcovers}
In order to formally define a near cover, we must first specify the larger set of maps that we allow. It is natural to restrict to surjective simplicial maps. We further restrict ourselves to $t$-to-$1$ maps and furthermore to the case where every edge is covered ``properly'' namely by a matching with exactly $t$ edges. This later restriction is not as arbitrary as it might seem because one can always modify a given map to have this property, without affecting any triangle that is properly covered. We denote the set of such maps by $M(X; t)$.

Formally, we will introduce a slightly more refined  definition. Let $G$ be a group acting on a set $S$ (such that $|S|=t$). We let $M(X; G,S)$ be the set of all pairs $(f',Y)$ such that
\begin{itemize}
\item $Y$ is a simplicial complex and $f':Y\to X$ is a surjective simplicial map.
\item For each $v\in X(0)$, $f^{-1}(v)$ can be identified with $S$.
\item For every edge $\set{u,v}\in X(1)$, there is a group element $g_{uv} \in G$, such that $f^{-1}(\set{u,v})$ is a bipartite matching between $f^{-1}(u)$ and $f^{-1}(v)$ viewed as two copies of $S$. This matching corresponds to the action of $g_{uv}$ on $S$. Namely, for every edge $\set{\tilde u, \tilde v}\in Y(1)$ such that $f(\set{\tilde u ,\tilde v})=\set{u,v}$ we have $g_{uv}(\tilde v) = (\tilde u)$.
\end{itemize}

We denote by $M_0(X;G,S) \subset M(X;G,S)$ the set of maps $(f,Y)$ that are genuine $(G,S)$-covers. By definition, this is the set of pairs with zero deficiency,
\[ M_0(X;G,S) = \Set{(f,Y)\in M(X;G,S)\; : \; m_f(Y)=0}.
\]

An important special case is when $S$ is a set of $t$ elements and $G$ the symmetric group acting on $S$, i.e. $G=Sym(S)$. In this case $M$ is simply $M(X;t)$ defined above, and $M_0$ becomes the set of all possible $t$-to-$1$ covers (with no restriction on the permutations covering any edge).
We will use shorthand $M$ and $M_0$ when the context is clear.

Inside $M$ we measure distance between two maps $(f_1,Y_1)$ and $(f_2,Y_2)$ by the fraction of edges $uv\in X(1)$ for which $f_1^{-1}(uv) \neq f_2^{-1}(uv)$. Note that comparing these two bipartite matchings makes sense through the natural identification  $Y_1(0) \longleftrightarrow (X(0)\times S) \longleftrightarrow Y_2(0) $.
This is a natural measure of distance as initiated in \cite{GGR} for testing of graph properties. In the context of two-dimensinoal complexes one could also compare the number of triangles that differ between $Y_1,Y_2$. However, the two distances are comparabe in our context because the weight of an edge is proportional to the number of triangles containing it, and the edge structure determines the allowed triangles for any map in $M$, so we focus on the edges:
\[ \dist((f_1,Y_1),(f_2,Y_2)) = \sum_{\{uv\in X(1) \;:\; f_1^{-1}(uv) \neq f_2^{-1}(uv)\}} c_X(uv) .
\]

We will be interested in the distance of $(f,Y)$ from being a genuine cover,
\[\dist((f,Y),M_0) = \min_{(f',Y')\in M_0}\, \dist((f,Y),(f'Y')).\]
We define the $(G,S)$-cover-stability to be the minimal ratio between the deficiency of $(f,Y)$ and its distance to a genuine cover. Let
\[ c(f,Y) = \frac{m_f(Y)}{\dist((f,Y),M_0)} .\]
The $GS$-cover-stability of $X$ is defined as
\[
c(X; G,S) = \min_{(f,Y)\in M\setminus M_0} c(f,Y)
\]
where of course both $M$ and $M_0$ here are taken with respect to $G$ and $S$.

\subsection{Non-Abelian First Cohomology}
\label{subs:nonab}

Let $X$ be a finite simplicial complex and let $G$ be a multiplicative group.
Let $C^0(X;G)$ denote the group of  $G$-valued functions on $X(0)$ with pointwise multiplication,
and let $$C^1(X;G)=\{\phi:X_{\ord}(1) \rightarrow G: \phi(u,v)=\phi(v,u)^{-1}\}.$$
The $0$-coboundary operator $d_0:C^0(X;G) \rightarrow C^1(X;G)$ be given by
$$d_0 \psi(u,v)=\psi(u)\psi(v)^{-1}.$$
For $\phi \in C^1(X;G)$ and $(u,v,w) \in X(2)$ let
$$d_1 \phi (u,v,w)=\phi(u,v)\phi(v,w)\phi(w,u).$$
Note that if $d_1\phi(u_1,u_2,u_3)=1$, then $d_1\phi(u_{\pi(1)},u_{\pi(2)},u_{\pi(3)})=1$ for all permutations $\pi$.
The set of \emph{$G$-valued $1$-cocycles} of $X$ is given by
$$Z^1(X;G)=\{\phi \in C^1(X;G):d_1\phi(u,v,w)=1 {\rm ~for~all~}(u,v,w) \in X_{\ord}(2)\}.$$
Define an action of $C^0(X;G)$ on $C^1(X;G)$ as follows.
For $\psi \in C^0(X;G)$ and  $\phi \in C^1(X;G)$ let
$$\psi . \phi (u,v)=\psi(u) \phi(u,v) \psi(v)^{-1}.$$
Note that $d_0 \psi =\psi . 1$ and that $Z^1(X;G)$ is invariant under the action of $C^0(X;G)$.
For $\phi \in C^1(X;G)$ let $[\phi]$ denote the orbit of $\phi$ under the action of $C^0(X;G)$.
The \emph{first cohomology of $X$ with coefficients in $G$} is the set of orbits
$$H^1(X;G)=\{[\phi]: \phi \in Z^1(X;G) \}.$$

\subsection{Correspondence of $1$-Cocycles and Covering Maps}
\label{subs:ccm}

We next recall the following classical construction (See Steenrod \cite{Steenrod51} for general spaces, and Surowski \cite{Surowski84} for the simplicial version).
Suppose $G$ acts on the left on a finite set $S$. For a $1$-cochain $\phi \in C^1(X;G)$, let $Y_{\phi}$ be the simplicial complex on the vertex set $Y_{\phi}(0)=\{[u,s]: u \in X(0), s \in S\}$, whose $k$-simplices are
$\tau=\big\{[u_0,s_0],\ldots,[u_k,s_k]\big\}$, where $\{u_0,\ldots,u_k\} \in X(k)$, and
$s_i=\phi(u_i,u_j)s_j$ for all $0 \leq i,j \leq k$. Let
$f:Y_{\phi} \rightarrow X$ be the simplicial projection map given by $f([u,s])=u$.
Note that if $\psi \in C^0(X;G)$, then there is an isomorphism $Y_{\psi . \phi} \cong_X Y_{\phi}$ via the simplicial map
$[v,s] \rightarrow [v,\psi(v)^{-1}s)]$.
\begin{theorem}[correspondence of cocycles and covers \cite{Surowski84}]
Let $X$ be a connected complex.
If $\phi \in Z^1(X;G)$ then $f:Y_{\phi} \rightarrow X$ is a covering map.
Conversely, let $f:Y \rightarrow X$ be a simplicial covering map and let $v_0 \in X(0)$. Then there is an action of $G=\pi_1(X,v_0)$ on $S=f^{-1}(v_0)$, and a $\phi \in Z^1(X;G)$ such that
$Y \cong_X Y_{\phi}$.
\end{theorem}

This correspondence extends to a correspondence between cochains and maps in $M(X; G,S)$. Indeed a map  $(f,Y)\in M(X;G,S)$ corresponds to a cochain $\phi \in C^1(X,G)$ such that for each edge $uv\in X(1)$ we have $\phi(uv)=g_{uv}$, where $g_{uv}$ is the group element that corresponds to the matching in $Y$ between $f^{-1}(u)$ and $f^{-1}(v)$.


\subsection{Cosystolic $1$-Expansion}
\label{subs:csy}

For a cochain $\phi \in C^1(X;G)$ let
$$
\supp(\phi)=\big\{\{u,v\} \in X(1): \phi(u,v) \neq 1\big\}
$$
and
$$
\supp(d_1\phi)=\big\{ \{u,v,w\} \in X(2) : d_1\phi(u,v,w)\neq 1 \big\}.
$$
Let
$$\|\phi\|=\sum_{e \in \supp(\phi)} c_X(e)$$
and
$$
\|d_1\phi\|=\sum_{\sigma \in \supp(d_1 \phi)} c_X(\sigma).$$
This measures the measure of triangles $uvw$ on which $d_1\phi(uvw) = \phi(uv)+\phi(vw)+\phi(wu)\neq 1$. For such a triangle we sometimes say that its equation isn't satisfied by $\phi$.
The \emph{distance} between $\phi, \psi \in C^1(X;G)$ is the measure of edges on which $\phi(e)\neq \psi(e)$,
$$
\dist(\phi,\psi)=\|\phi\psi^{-1}\|.
$$
The \emph{cosystolic norm} of $\phi \in C^1(X;G)$ is the distance of $\phi$ from $Z^1(X;G)$, i.e.
$$\|\phi\|_{\csy}=\min \{\|\phi \psi^{-1}\|: \psi \in Z^1(X;G)\}.$$
This is measuring the distance of $\phi$ to the closest cocycle $\psi\in Z^1(X;G)$, in terms of how many edges need to be changed to go from $\phi$ to $\psi$, and taking into account the weights of the edges.
The \emph{cosystolic expansion} of $\phi \in C^1(X;G)\setminus Z^1(X;G)$ is
$$
h(\phi)=\frac{\|d_1\phi\|}{\|\phi\|_{\csy}}.
$$
The \emph{cosystolic expansion} of $X$ is
$$
h_1(X;G)=\min\big\{h(\phi): \phi \in C^1(X;G) \setminus Z^1(X;G)\big\}.
$$
When this is at least a constant, it means that if $\|d_1\phi\|$ is small, namely $\phi$ satisfies most of the triangle equations, then $\phi$ is at most $\frac{\|d_1\phi\|}{h(X;G)} = O(\|d_1(\phi)\|)$-close to a genuine cocycle.

\noindent
{\bf Example:} Let $\dn$ denote the $(n-1)$-simplex. In \cite{M13} it is shown that for any group $G$
$$
h_1(\dn;G) \geq \frac{n}{n-2} > 1.
$$

 \subsection{Cover-stability and cosystolic expansion}
\label{subs:res}

In this section we are finally ready to formally state our results.
Our main result is that a simplicial complex $X$ is a cosystolic expander with respect to $G$ if and only if $X$ is $(G,S)$-cover stable.

We first need one more definition. For $g \in G$,
let  $\fix(g)=|\{s \in S: gs=s\}|$.
The \emph{fixity} of the action of $G$ on $S$ is $\ffix_G(S)= \max_{g \neq 1} \fix(g)$.
The action of $G$ is \emph{faithful} if $\ffix_G(S) < |S|$, and in this case clearly $\ffix_G(S)\leq |S|-2$. The action of $G$ is \emph{free} if $\ffix_G(S)=0$.
\begin{theorem}[stability $\leftrightarrow$ expansion]
\label{t:mainformal}
Let $X$ be a simplicial complex. Let $G$ act on a finite set $S$.
\[\frac 2 {|S|}\cdot h_1(X;G) \le \big(1-\frac{\ffix_G(S)}{|S|}\big) \cdot h_1(X;G) \le c(X;G,S) \le h_1(X;G)
\]
In particular, $X$ is $(G,S)$-cover-stable iff it is a cosystolic expander with respect to $G$ coefficients.
\end{theorem}
This theorem follows immediately from the following.
Let $\phi \in C^1(X;G)$.
The following result shows, roughly speaking, that if the deficiency of $f:Y_{\phi} \rightarrow X$ is small,
then $\phi$ is close to a $1$-cocycle in $H^1(X;G)$ and therefore $Y_\phi$ is close to a genuine cover.
\begin{theorem}
\label{t:nearc}
Let $G$ act on a finite set $S$.
Then for any $\phi \in C^1(X;G)$ there exists a $\psi \in Z^1(X;G)$ such that
\begin{equation}
\label{e:nearc}
\dist(\phi,\psi)
 \leq \frac{m(Y_{\phi})}{\big(1-\frac{\ffix_G(S)}{|S|}\big) \cdot h_1(X;G)} .
\end{equation}
\end{theorem}

\subsection{An example for a cover-stable complex: the spherical building} \label{subs:sph}
Let $X$ be the spherical building $A_{3}(\FF_q)$, i.e. the order complex of the lattice of all nontrivial
linear subspaces of $\FF_q^4$.
\begin{theorem}
\label{t:sph}
For any finite group $G$
$$
h_1\big(A_3(\FF_q);G \big) \geq \frac{1}{9}.
$$
\end{theorem}
\noindent
Combining Theorems \ref{t:nearc} and \ref{t:sph} we obtain the following
\begin{corollary}
\label{c:sph}
Let $G$ act on a finite set $S$. Then for
any $\phi \in C^1(A_3(\FF_q);G)$ there exists a $\psi \in Z^1(A_3(\FF_q);G)$ such that
$$\dist(\phi,\psi) \leq \frac{9 \, m(Y_{\phi})}{\big(1-\frac{\ffix_G(S)}{|S|}\big)}.$$
In particular, if the action of $G$ is free then
$$\dist(\phi,\psi) \leq 9\, m(Y_{\phi}).$$
\end{corollary}
\noindent
{\bf Remark:} The simple connectivity of $A_3(\FF_q)$ implies that all cocycles of this complex are coboundaries, i.e. the cohomology vanishes. This means that there are no non-trivial covers, so if $\psi \in Z^1(X;G)$, then $Y_{\psi}$ is isomorphic to the trivial $|S|$-fold covering of  $A_3(\FF_q)$.
\ \\ \\
The remaining of this note is organized as follows. In Section \ref{s:nearc} we prove our main results, Theorems \ref{t:mainformal} and \ref{t:nearc}.
In Section \ref{s:egl} we obtain Theorem \ref{t:sph}, as a consequence of a general bound (Theorem \ref{t:latticelb}) on the non-abelian $1$-expansion of order complexes of geometric lattices.

\section{Proof of Theorems \ref{t:mainformal} and \ref{t:nearc}}
\label{s:nearc}
\paragraph{Proof of Theorem  \ref{t:nearc}.}
Let $\phi \in C^1(X;G) \setminus Z^1(X;G)$. Recall that $f:Y_{\phi} \rightarrow X$ is the projection map $f([u,s])=u$. Let $u \in X(0)$ and let $e=\{v_1,v_2\} \in \lk(X,u)$. Then
$e=\{v_1,v_2\} \in D_f([u,s])$ iff
$$
\{[u,s],[v_1,\phi(v_1,u)s],[v_2,\phi(v_2,u)s]\}  \not\in Y_{\phi}(2),
$$
i.e. iff
\begin{equation}
\label{e:dphi}
d_1\phi(u,v_1,v_2)s \neq s.
\end{equation}
Using (\ref{e:dphi}) and (\ref{e:prod1}) we obtain
\begin{equation}\label{eq:chain}
\begin{split}
|S|\,m(Y_{\phi})&=\sum_{u \in X(0)}c_X(u) \sum_{\tilde{u} \in f^{-1}(u)} \mu_f(\tilde{u}) \\
&=\sum_{u \in X(0)}\sum_{s \in S} c_X(u) \mu([u,s]) \\
&=\sum_{u \in X(0)} c_X(u) \sum_{s \in S} \sum_{e \in D_f([u,s])} c_{\lk(X,u)}(e) \\
&=\sum_{u \in X(0)} c_X(u) \sum_{\{v_1,v_2\} \in \lk(X,u)} |\{s:d_1\phi(u,v_1,v_2)s \neq s\}|  \cdot c_{\lk(X,u)}(\{v_1,v_2\}) \\
&=\sum_{u \in X(0)} c_X(u) \sum_{e \in \lk(X,u)(1)} \big(|S|-\fix\left(d_1\phi(u \cup e)\right)\big) \cdot c_{\lk(X,u)}(e)
\end{split}
\end{equation}
At this point we bound the expression from above and from below. For the lower bound,
\begin{equation*}
\begin{split}
\eqref{eq:chain}&\geq \big(|S|-\ffix_G(S)\big)\sum_{u \in X(0)}
\sum_{\{e \in \lk(X,u)(1):d_1\phi(u \cup e) \neq 1\}}
c_X(u)c_{\lk(X,u)}(e) \\
&= \big(|S|-\ffix_G(S)\big)\sum_{u \in X(0)}
\sum_{\{e \in \lk(X,u)(1):u \cup e \in \supp(d_1\phi)\}}
\frac{1}{3} \cdot c_X(u \cup e)  \\
&=\big(|S|-\ffix_G(S)\big)\|d_1\phi\|.
\end{split}
\end{equation*}
Next, for the upper bound,
\begin{equation*}
\begin{split}
\eqref{eq:chain}&\leq |S|\,\sum_{u \in X(0)}
\sum_{\{e \in \lk(X,u)(1):d_1\phi(u \cup e) \neq 1\}}
c_X(u)c_{\lk(X,u)}(e) \\
&= |S|\,\sum_{u \in X(0)}
\sum_{\{e \in \lk(X,u)(1):u \cup e \in \supp(d_1\phi)\}}
\frac{1}{3} \cdot c_X(u \cup e)  \\
&=|S|\cdot\|d_1\phi\|.
\end{split}
\end{equation*}
We conclude that
\begin{equation}\label{eq:ineq} \frac{|S|-\ffix_G(S)}{|S|}\cdot  \|d_1\phi\| \le m(Y_\phi) \le \|d_1\phi\|.
\end{equation}
\noindent
As $h_1(X;G) \leq \frac{\|d_1\phi\|}{\|\phi\|_{\csy}}$, it follows that
\begin{equation*}
\label{e:agp}
\begin{split}
\min &\{\dist(\phi,\psi):\psi \in Z^1(X;G)\}= \|\phi\|_{\csy} \\
&\leq \frac{\|d_1\phi\|}{h_1(X;G)} \leq \frac{|S|m(Y_{\phi})} {\big(|S|-\ffix_G(S)\big)h_1(X;G)} \\
&=\frac{m(Y_{\phi})}{\big(1-\frac{\ffix_G(S)}{|S|}\big) \cdot h_1(X;G)} ~~.
\end{split}
\end{equation*}
{\enp}

\paragraph{Proof of Theorem  \ref{t:mainformal}.}
We deduce this theorem from the proof of Theorem \ref{t:nearc}. The first inequality is trivial because a group element acting non-trivially on $S$ must have $Fix_G(S)\leq |S|-2$. The remaining inequalities follow from the correspondence between cochains $\phi \in C^1(X;G)$ and maps $(f,Y_\phi)\in M(X;G,S)$. Recall that
\[ h(X;G,S) = \min_\phi \frac{\|d_1\phi\|}{\dist(\phi,Z^1(X;G))}, \qquad
c(X;G,S) = \min_{(f,Y)} \frac{m_f(Y)}{\dist((f,Y),M_0)}.
\]
For every cochain $\phi$ the denominators are the identical: $\dist((f,Y_\phi), M_0) = \dist(\phi,Z^1(X;G))$. The numerators satisfy \eqref{eq:ineq}. Taking minimum over all $\phi$ completes the proof.
\qed

\section{The $1$-Expansion of Geometric Lattices}
\label{s:egl}

Let $(P,\leq)$ be a finite poset. The \emph{order complex} of $P$ is
the simplicial complex on the vertex set $P$ whose simplices are the
chains $a_0<\cdots < a_k$ of $P$. In the sequel we identify a poset
with its order complex.
A poset $(L,\leq)$ is a \emph{lattice} if any two elements $x,y \in L$ have a unique minimal upper bound $x \vee y$ and a unique maximal lower bound $x \wedge y$. A lattice $L$ with minimal element $\hzero$ and maximal element $\hone$ is \emph{ranked}, with rank function
$\rk(\cdot)$, if $\rk(\hzero)=0$ and $\rk(y)=\rk(x)+1$ whenever $y$ covers $x$. $L$ is a \emph{geometric lattice} if $\rk(x)+\rk(y) \geq \rk(x \vee y)+\rk(x \wedge y)$ for any $x,y \in L$,
and any element in $L$ is a join of atoms (i.e., rank $1$ elements).

\begin{figure}
\begin{center}
  \scalebox{0.4}{\input{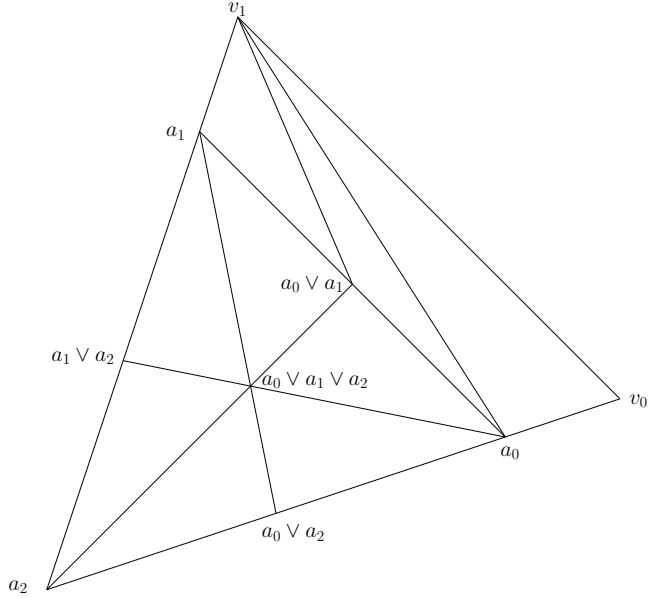}}
  \caption{The subcomplex $Y_s(uv)$.}
  \label{figure1}
\end{center}
\end{figure}

Let $L$ be a geometric lattice with $\rk(\hone)=n \geq 3$. A classical result of Folkman \cite{F66} asserts that $\overline{L}=L-\{\hzero,\hone\}$ is homotopy equivalent to a wedge of $(n-2)$-spheres. In particular, $\overline{L}$ is simply connected, and hence $H^1(\ovl;G)=\{1\}$ for any group $G$. Here we provide a lower bound for $h_1(\ovl;G)$.
Let $S$ be a set of linear orderings on the set of atoms $A$ of $L$, equipped with a probability distribution $\mu$.
Let $\prec_{s}$ denote the ordering associated with $s \in S$. For $s \in S$ and $u \in \ovl$ let
$a(s)=\min A$, $b(s,v)=\min\{a \in A:a \leq u\}$ where both minima are taken with respect to
$\prec_s$. For $s \in S$ and $v_0<v_1 \in \ovl$, let $a_0=b(s,v_0), a_1=b(s,v_1), a_2=a(s)$. Clearly
$a_2 \preceq_s a_1 \preceq_s a_0$. Let $Y_s(v_0v_1)$ be the $2$-dimensional subcomplex of $\ovl$ depicted in Figure \ref{figure1}.
For $\tau \in \ovl(2)$ and $s \in S$ let
\begin{equation}
\label{e:dst}
\delta_s(\tau):=\sum_{\{uv \in \ovl(1): \tau \in Y_s(uv)\}} \frac{c_X(uv)}{c_X(\tau)}.
\end{equation}
Let $\delta(\tau)=E[\delta_s(\tau)]$ denote the expectation of $\delta_s(\tau)$.
The proof of the next result uses a homotopical adaptation of the approach of \cite{KM18}.
\begin{theorem}
\label{t:latticelb}
\begin{equation*}
\label{e:latticelb}
h_1(\ovl;G) \geq \left( \max_{\tau \in \ovl(2)} \delta(\tau) \right)^{-1}.
\end{equation*}
\end{theorem}
\noindent
We will need the following simple fact.
\begin{claim}
\label{c:simc}
Let $G$ be a group and
let $K$ be a $2$-dimensional simply connected simplicial complex. Suppose $v_0,\ldots,v_{m-1},v_m=v_0$ are the vertices of a $1$-cycle in $K$. If $\phi \in C^1(K;G)$ satisfies $$\phi(v_0,v_1)\cdot\phi(v_1,v_2) \cdots \phi(v_{m-1},v_0) \neq 1$$ then there exists a $2$-simplex $(a,b,c) \in K_{\ord}(2)$ such that $d_1\phi(a,b,c) \neq 1$.
{\enp}
\end{claim}
\noindent
{\bf Proof of Theorem \ref{t:latticelb}:}
Let $\phi \in C^1(\ovl;G)$. For $s \in S$, define $\psi_s \in C^0(\ovl;G)$ by
$$
\psi_s(u)=\phi\big(a(s),a(s)\vee b(s,u)\big)\cdot\phi\big(a(s) \vee b(s,u),b(s,u) \big)\cdot \phi\big(b(s,u),u\big).$$
Let $v_0<v_1 \in \ovl$ and (as before) denote $a_0=b(s,v_0), a_1=b(s,v_1), a_2=a(s)$.
Consider the $1$-cycle in $Y_s(uv)$ whose vertices are
$$
(x_0,\ldots,x_7)=(a_2,a_0 \vee a_2,a_0,v_0,v_1,a_1,a_1 \vee  a_2,a_2).
$$
Then
\begin{equation*}
\label{e:prcyc}
\begin{split}
(\psi_s).\phi(v_0,v_1)&=\psi_s(v_0)\phi (v_0,v_1) \psi_s(v_1)^{-1} \\
&=\phi(x_0,x_1)\phi(x_1,x_2)\cdots \phi(x_5,x_6)\phi(x_6,x_0).
\end{split}
\end{equation*}
\noindent
Since $Y_s(uv)$ is simply connected (in fact contractible), it follows from Claim \ref{c:simc} that if
$(\psi_s).\phi(v_0,v_1) \neq 1$, then there exists a $2$-simplex $(x,y,z) \in Y_s(uv)$ such that $d_1\phi(x,y,z) \neq 1$.
Therefore
\begin{equation}
\label{e:bndh}
\begin{split}
\|\phi\|_{\csy} &\leq \sum_{s \in S} \mu(s) \|(\psi_s).\phi\| \\
&=\sum_{s \in S} \mu(s) \sum \{c_X(uv): uv \in \ovl(1), (\psi_s).\phi(u,v) \neq 1\} \\
&\leq \sum_{s \in S} \mu(s) \sum \{c_X(uv): uv \in \ovl(1), \supp(d_1\phi) \cap Y_s(uv) \neq \emptyset\} \\
&\leq \sum_{s \in S} \mu(s) \sum_{\tau \in \supp(d_1\phi)}
\sum \{c_X(uv): uv \in \ovl(1), \tau \in Y_s(uv)\} \\
&=\sum_{\tau \in \supp(d_1\phi)} c_X(\tau) \sum_{s \in S} \mu(s)
\sum \{\frac{c_X(uv)}{c_X(\tau)}: uv \in \ovl(1), \tau \in Y_s(uv)\} \\
&= \sum_{\tau \in \supp(d_1\phi)} c_X(\tau) E[\delta_s(\tau)] \\
&=\sum_{\tau \in \supp(d_1\phi)} c_X(\tau) \delta(\tau)
\leq \|d_1\phi\| \max_{\tau \in \ovl(2)} \delta(\tau).
\end{split}
\end{equation}
{\enp}
\noindent
For lattices $L$ with sufficient symmetry (e.g. spherical buildings), Theorem \ref{t:latticelb} can be used to give explicit lower bounds on $h_1(\ovl,G)$.
\ \\ \\
{\bf Proof of Theorem \ref{t:sph}:}
Let $L$ be the lattice of all nontrivial linear subspaces of $\FF_q^4$. Then $\ovl= A_3(\FF_q)$.
Let $\prec$ be an arbitrary fixed linear order on the set of atoms $A$.
Let $S$ be the group $GL_4(\FF_q)$ with the uniform distribution. For $s \in S$ let $\prec_s$ be the linear order
on $A$ given by $a \prec_s a'$ if $s^{-1}a \prec s^{-1}a'$.
Let $\idd$ denote the identity element of $S$. It is straightforward to check that
$a(\idd)=s^{-1}a(s)$ and $b(\idd,u)=s^{-1}b(s,su)$ for any $u \in \ovl(0)$.
Hence $Y_s(e)=sY_{\idd}(s^{-1}(e))$ for any $e \in \ovl(1)$. It follows that if
$e \in \ovl(1)$, $\tau \in \ovl(2)$ and $t \in S$, then
$t(\tau) \in Y_s(e)$ iff $\tau \in Y_{t^{-1}s}(t(e))$.
This implies that for any $t \in S$ and $\tau \in \ovl(2)$
\begin{equation}
\label{e:equsets}
|\{(s,e)\in S\times \ovl(1): \tau \in Y_s(e)\}|=
|\{(s,e)\in S\times \ovl(1): t(\tau) \in Y_s(e)\}|.
\end{equation}
Since $$\frac{c_{\ovl}(e)}{c_{\ovl}(\tau)}=\frac{q+1}{3},$$ it follows from (\ref{e:dst}) and (\ref{e:equsets}) that
$\delta(\tau)=\delta(t(\tau))$ for any $t \in S$ and $\tau \in \ovl(2)$. Next note that $S$ is transitive
on $\ovl(2)$ and thus $\delta$ is constant, i.e. $\delta(\tau)=\gamma$ for all $\tau \in \ovl(2)$. Therefore
\begin{equation*}
\begin{split}
f_2(\ovl) \gamma &= \sum_{\tau \in \ovl(2)} \delta(\tau) \\
&=\sum_{\tau \in \ovl(2)} \frac{1}{|S|} \sum_{s \in S}\delta_s(\tau) \\
&=\frac{q+1}{3|S|}|\{(s,uv,\tau) \in S \times \ovl(1) \times \ovl(2): \tau \in Y_s(uv)\}| \\
&=\frac{q+1}{3|S|} \sum_{s \in S} \sum_{uv \in \ovl(1)} f_2(Y_s(uv)) \\
&\leq \frac{q+1}{3|S|} \cdot |S| \cdot f_1(\ovl)\cdot 9=3(q+1) f_1(\ovl).
\end{split}
\end{equation*}
Therefore
$$
\gamma \leq \frac{3(q+1) f_1(\ovl)}{f_2(\ovl)}=9,
$$
hence Theorem \ref{t:sph} follows from Theorem \ref{t:latticelb}.
{\enp}


\begin{thebibliography}{99}


\bibitem{ChapmanLP19}
Michael Chapman, Nati Linial, and Yuval Peled.
\newblock Expander graphs -- both local and global.
\newblock In {\em 59th IEEE Annual Symposium on Foundations of Computer Science
  (FOCS)}. IEEE Canada, 2019.

\bibitem{Conlon2019}
David Conlon.
\newblock Hypergraph expanders from cayley graphs.
\newblock {\em Israel Journal of Mathematics}, Jul 2019.

\bibitem{dinur2017high}
Irit Dinur and Tali Kaufman.
\newblock High dimensional expanders imply agreement expanders.
\newblock In {\em 58th IEEE Annual Symposium on Foundations of Computer Science
  (FOCS)}. IEEE Canada, 2017.

\bibitem{DinurKamber}
Irit Dinur and Amitay Kamber.
\newblock Unpublished manuscript.

\bibitem{EvraK16}
Shai Evra and Tali Kaufman.
\newblock Bounded degree cosystolic expanders of every dimension.
\newblock In {\em Proceedings of the 48th Annual {ACM} {SIGACT} Symposium on
  Theory of Computing, {STOC} 2016, Cambridge, MA, USA, June 18-21, 2016},
  pages 36--48, 2016.

\bibitem{F66}
J.\ Folkman, The homology groups of
a~lattice, {\it J.\ Math.\ Mech.}, {\bf 15}(1966) 631--636.


\bibitem{Garland73}
H. Garland, $p$-adic curvature and the cohomology of discrete subgroups of p-adic groups,
{\it Ann. of Math.},{\bf 97}(1973) 375--423.

\bibitem{GGR}
Oded Goldreich, Shafi Goldwasser, and Dana Ron, Property Testing and its Connection to Learning and Approximation, {\it J. ACM} {\bf 45(4)}(1998) 653--750.

\bibitem{Gromov10}
M. Gromov, Singularities, expanders and topology of maps. Part 2: From combinatorics to topology via algebraic isoperimetry, {\it Geom. Funct. Anal.} {\bf 20}(2010) 416-�526.


\bibitem{KKL14}
Tali Kaufman, David Kazhdan, and Alexander Lubotzky.
\newblock Ramanujan complexes and bounded degree topological expanders.
\newblock In {\em 55th {IEEE} Annual Symposium on Foundations of Computer
  Science, {FOCS} 2014, Philadelphia, PA, USA, October 18-21, 2014}, pages
  484--493, 2014.

\bibitem{KaufmanL14}
Tali Kaufman and Alexander Lubotzky.
\newblock High dimensional expanders and property testing.
\newblock In {\em Innovations in Theoretical Computer Science, ITCS'14,
  Princeton, NJ, USA, January 12-14, 2014}, pages 501--506, 2014.

\bibitem{KM18}
D. N. Kozlov and R. Meshulam,
Quantitative aspects of acyclicity, {\it Research in the Mathematical Sciences}, to appear.


\bibitem{LM06} N. Linial and R. Meshulam, Homological connectivity
of random $2$-complexes, {\it Combinatorica} {\bf 26}(2006)
475--487.

\bibitem{siqi}
Siqi Liu, Sidhanth Mohanty, and Elizabeth Yang.
\newblock High-dimensional expanders from expanders.
\newblock {\em CoRR}, abs/1907.10771, 2019.

\bibitem{M13}
R. Meshulam, Bounded quotients of the fundamental group of a random 2-complex,
arXiv:1308.3769 .


\bibitem{MW09}
R. Meshulam and N. Wallach, Homological connectivity of
random k-dimensional complexes, {\it Random Structures and Algorithms} {\bf 34}(2009) 408--417.



\bibitem{Olum58}
P. Olum, Non-abelian cohomology and van Kampen's theorem, {\it Ann. of Math.}, {\bf 68}(1958) 658�-668.


\bibitem{OppGarland}
Izhar Oppenheim.
\newblock Local spectral expansion approach to high dimensional expanders part
  {I:} descent of spectral gaps.
\newblock {\em Discrete {\&} Computational Geometry}, 59(2):293--330, 2018.


\bibitem{Steenrod51}
N. Steenrod,
The Topology of Fibre Bundles. Princeton Mathematical Series, vol. 14. Princeton University Press, Princeton, N. J., 1951.

\bibitem{Surowski84}
D. B. Surowski, Covers of simplicial complexes and applications to geometry, {\it Geom. Dedicata},
{\bf 16}(1984) 35-�62.

\bibitem{RVW}
Omer Reingold, Salil Vadhan, and Avi Wigderson.
\newblock Entropy waves, the zig-zag graph product, and new constant-degree
  expanders and extractors.
\newblock {\em Annals of Mathematics}, 155:157--187, 2002.


\end{thebibliography}

\end{document}